\documentclass[final]{ifacconf}

\usepackage[utf8x]{inputenc}
\usepackage{amssymb,amsmath}
\usepackage{mathtools}
\usepackage{mathrsfs}
\usepackage{amssymb}
\usepackage{graphicx}          % Include this line if your
\theoremstyle{plain}
\newtheorem{theorem}{Theorem}
\newtheorem{assumption}[theorem]{Assumption}
\newtheorem{lemma}[theorem]{Lemma}
\newtheorem{definition}[theorem]{Definition}
\newtheorem{proposition}[theorem]{Proposition}

\usepackage[usenames,dvipsnames]{xcolor}
\usepackage{soul}

\renewcommand{\Re}{\mathbb{R}}
\newcommand{\uu}{\mathbf{u}}
\newcommand{\N}{\mathbb{N}}
\newcommand{\NN}{\mathcal{N}}
\newcommand{\FF}{\mathfrak{F}}
\newcommand{\expect}{\mathbb{E}}
\renewcommand{\prob}{\mathrm{P}}

\newcommand{\EE}{\mathscr{E}}

\newcommand{\dfn}{\mathrel{:}=}
\newcommand{\bet}{\mathrm{bet}}

\renewcommand{\Re}{\mathbb{R}}

\begin{document}

\begin{frontmatter}
\title{Stochastic economic model predictive control for Markovian switching systems%
\thanksref{footnoteinfo}} % Title, preferably not more than 10 words.

%
% ACKNOWLEDGEMENT OF FUNDING
%
\thanks[footnoteinfo]{This work was supported by the EU-funded H2020
research project DISIRE, grant agreement No.~636834.}

\author[KUL]{Pantelis Sopasakis},
\author[IMTLUCCA]{Domagoj Herceg},
\author[KUL]{Panagiotis Patrinos} and
\author[IMTLUCCA]{Alberto Bemporad}

\address[IMTLUCCA]{IMT School for Advanced Studies Lucca, Piazza San Ponziano 6, 55100 Lucca, Italy.}
\address[KUL]{KU Leuven, Dept. of Electrical Engineering (ESAT-STADIUS), Kasteelpark Arenberg 10, 3001 Leuven, Belgium.}

\begin{abstract}  % Abstract of not more than 250 words.
The optimization of process economics within the model predictive control (MPC)
formulation has given rise to a new control paradigm known as economic MPC (EMPC). 
Several authors have discussed the closed-loop properties of EMPC-controlled deterministic systems, 
however, little have uncertain systems been studied. In this paper we propose EMPC 
formulations for nonlinear Markovian switching systems which guarantee recursive feasibility, 
asymptotic performance bounds and constrained mean square (MS) stability.
\end{abstract}

\end{frontmatter}
 
\section{Introduction}
\subsection{Background and motivation}
Recently, a new approach to model predictive control (MPC) termed 
\textit{economic model predictive control} (EMPC) has gained a lot 
of attention. Rather than minimizing a deviation from a prescribed 
(optimal/best) set-point or a tracking reference, the main objective 
in EMPC is to optimize a given economic cost functional~\citep{AngAmrRaw_TAC}. 
Often, in engineering practice, the main objective is to devise control
algorithms which asymptotically guarantee an economic operation of the 
controlled plant. 

Already, a considerable body of theoretical results has been
reported in the literature characterizing the asymptotic performance 
of EMPC. Perhaps \textit{dissipativity} is the most salient notion 
in the pertinent literature which is shown to be a sufficient condition for 
proving optimal operation at a steady state and stability of EMPC formulations%
~\citep{AngAmrRaw_TAC}. The same authors show that economic MPC has no 
worse an asymptotic average performance than the best admissible steady 
state operation --- however, the converse is not true~\citep{muller2013convergence}.

The introduction of a, possibly non-quadratic and nonconvex, 
economic cost into the MPC framework disqualifies the standard 
stability analysis used in the MPC literature.
\cite{AngAmrRaw_TAC} propose the use of a simple terminal constraint
to guarantee stability of EMPC-controlled systems which is generalized
by~\cite{AmrRawAng_ARC} using terminal set constraints. 
\cite{fagiano-generalized} use a generalized terminal state constraint, 
where terminal state constraint is left as a free variable
to be optimized which helps to increase feasibility region of EMPC. 
This concept was further generalized to include terminal region 
constraint~\citep{muller2014performance}.
It was further shown that EMPC can achieve near-optimal operation without terminal
constraints and costs for a sufficiently large prediction horizon~\citep{grune}.
Similar results exist for a system that is best operated 
at a periodic regime~\citep{zanon2013}. It is worth noting that 
this wealth of results concerns only deterministic systems.

In spite of the noticeable interest for the idea of EMPC 
there are very few theoretical results accounting for uncertainty 
which is inevitable in a real-world operation.
% @Pantelis: These guys do nothing even remotely related to EMPC...
% [...] ~\citep{jorgensen-power-systems}.
\cite{Bo20149412} propose a scenario-based EMPC formulation for 
fault-tolerant constrained regulation and a similar approach is pursued
by~\cite{luc+14}.
\cite{Lucia_2014}  present a multi-stage scenario-based nonlinear MPC 
control strategy validated on a benchmark example,
but no performance guarantees or stability analysis is provided.
An interesting theoretical treatment is given by~\cite{bayer-JPC-2014}
where a tube-based EMPC formulation is proposed for constrained systems 
with bounded additive disturbances. 
Very recently \cite{Bayer2016151} proposed a robust economic MPC formulation 
for linear systems with bounded additive uncertainty with known probability 
distribution.

\subsection{Contributions} 
In this paper we endeavor to cover the theoretical gap in EMPC
for an important class of stochastic systems --- the Markovian switching systems.
We first study the properties of an MPC formulation for Markovian switching systems
where optimal steady states are mode-dependent. We propose an MPC
scheme which is recursively feasible and satisfies an asymptotic performance bound.
Assuming that there is a common optimal steady state, we show that the MPC-controlled
system is mean-square (MS) stable when a stochastic dissipativity condition is satisfied.
We then formulate a variant of the MPC problem using mode-dependent terminal 
constraints and provide mean-square stability conditions and performance bounds.
We then provide guidelines for the design of mean-square stabilizing predictive 
controllers for nonlinear systems imposing weak conditions on the system 
dynamics and the EMPC stage cost.

\subsection{Notation and mathematical preliminaries}
% INTRODUCTORY: R, R_+
Let $\Re$ and $\Re_+$, $\Re^n$, $\Re^{n\times n}$ denote 
the sets of real numbers, nonnegative reals, 
$n$-dimensional real vectors and $n$-by-$m$ matrices. 
Let $\mathcal{B}_\delta$ be the ball of $\Re$ of radius $\delta$, 
that is $\mathcal{B}_\delta \dfn \{x: \|x\| < \delta\}$.
% LSC
A function $f:\Re^n\to \Re$ is called \textit{lower semicontinuous} if
its epigraph, that is the set $\operatorname{epi}f = \{(x,\alpha)\in\Re^{n+1}: f(x) 
\leq \alpha\}$, is closed. 
% LEVEL BOUNDEDNESS
We say that $f:\Re^n\to \Re$ is 
\textit{level-bounded} if its level sets, $\operatorname{lev}_\alpha f = 
\{x: f(x) \leq \alpha\}$, are bounded. 
We say that $f:\Re^n\times \Re^m\ni (x,u) \mapsto f(x,u) \in \Re$ is 
% UNIFORM LEVEL BOUNDEDNESS
\textit{level-bounded in $u$ locally uniformly in $x$} if for every 
$\bar{x}$ there is a neighborhood of $\bar{x}$, $V_{\bar{x}}\subseteq \Re^n$, 
so that $\{(x,u): x\in V_{\bar{x}}, f(x,u) \leq \alpha\}$ is bounded.
% BETA-SMOOTH
A function $f:\Re^n\to\Re^m$ is called \textit{$\beta$-smooth} if
it is differentiable with $\beta$-Lipschitz gradient, that is
$\| \nabla f(y) - \nabla f(x)  \| \le \beta  \| y - x \|$
for all $x,y\in\Re^n$; then, we have that
$\| f(y) - f(x) - \nabla f(x) (y - x) \| \le 
\tfrac{\beta}{2} \| y - x \|^2$.
% POSITIVE DEFINITENESS
We say that a function $f:\Re^n\to\Re$ is positive definite
around $x_0$ if $f(x_0)=0$ and $f(x)>0$ for $x\neq x_0$.
$A\succcurlyeq 0$ denotes that $A$ is a positive semidefinite matrix
and $A\succ 0$ means that $A$ is positive definite.
We denote the transpose of a matrix $A$ by $A^\top$.

\section{Stochastic Economic Model Predictive Control}

\subsection{System dynamics}
Consider the following Markovian switching system
\begin{align}
 x_{k+1} = f(x_k, u_k, \theta_k),
 \label{eq:system_evo}
\end{align}
driven by the random parameter $\theta_k$ which is a time-homogeneous
irreducible and aperiodic Markovian process with values in a 
finite set $\NN=\{1,\ldots, \nu\}$
with transition matrix $P=(p_{ij})\in\Re^{\nu\times \nu}$ and 
\textit{initial distribution} $v=(v_1,\ldots, v_\nu)$~\citep{mjls2005costa}.
We assume that at time $k$ we measure the
full state $x_k$ and the value of $\theta_k$.
Markov jump linear systems (MJLS) with additive disturbances are a special case 
of~\eqref{eq:system_evo} with $f(x,u,\theta) = A_\theta x + B_\theta u + w_\theta$.

Let  $\Omega\dfn \prod_{k\in\N} (\Re^n\times\Re^m\times\NN)$
and $\FF_k$ be the minimal $\sigma$-algebra over the Borel-measurable
rectangles of $\Omega$ with $k$-dimensional base and
$\FF$ be the minimal $\sigma$-algebra over all Borel-measurable
rectangles.
Define the filtered probability space  $(\Omega,\FF,\{\FF_k\}_{k\in\N},\prob)$
where $\prob$ is the unique product probability
measure according to~\citep[Th. 2.7.2]{ash1972real}
with $\prob(\theta_{0}=i_0,\theta_{1}=i_1,
\ldots,\theta_{k}=i_k)=v_{i_0}p_{i_0i_1}\cdots p_{i_{k-1}i_k}$ for any $i_0,i_1,\ldots,i_k\in\NN$
and $k\in\N$, where $\theta_{k}$ is an $\FF_k$-adapted random variable from
$\Omega$ to $\NN$.
We will use the notation $u\lhd \FF_k$ to denote that the random variable $u$
is $\FF_k$-measurable.

Let $\expect[\cdot]$ denote the expectation of a random variable with respect to $\prob$ and
$\expect[\cdot|\FF_k]$ the conditional expectation.
It can be shown~\citep{tejada2010nonlinear}  that the augmented state
$(x_{k},\theta_{k})$ contains all the probabilistic information relevant
to the evolution of the Markovian switching system for times $t>k$.

%
%
% ..... DEFINITION:
%       Cover of a node
%       Bet node        ......
%
\begin{definition}[Cover and bet node]
 For every node $i\in\NN$, the \textit{cover} of $i$ is the set
 $\mathcal{C}(i) = \{j\in\NN\mid p_{ij}>0\}$. The \textit{bet node} of an $i\in\NN$
 is a node $\bet(i)\in\mathcal{C}(i)$ with $p_{i\bet(i)} \geq p_{ij}$
 for all $j\in\mathcal{C}(i)$.
\end{definition}

A bet of a mode $\theta_k = i$ is one of the most likely successor
modes $\theta_{k+1}$.

System~\eqref{eq:system_evo} is subject to the following
joint state-input  constraints
\begin{align}\label{eq:constraints}
 (x_k, u_k) \in Y_{\theta_k}.
\end{align}

%
%
% ------ ASSUMPTION
%        Well-posedness   ------
%
Let $\ell:\Re^n\times\Re^m\times\NN\to \Re$ be a mode-dependent 
cost function. 
\begin{assum}[Well-posedness]\label{assum:well-posedness}
 % assumptions on ell
 For each $\theta\in\NN$, $\ell(\cdot, \cdot, \theta)$ are
 nonnegative, lower semicontinuous and level-bounded in $u$ 
 locally uniformly in $x$, 
 % assumptions on f
 $f(\cdot, \cdot, \theta)$ are continuous and 
 % assumptions on Y
 sets $Y_\theta$ are nonempty and compact.
 % assumptions on the random process
 The random process $\{\theta_k\}_k$ is an irreducible and
 aperiodic Markov chain.
\end{assum}

%
%
% ------ ASSUMPTION
%        Optimal steady states   ------
%
\begin{definition}[Optimal steady states]
\label{def:optimal-steady-states}
Given a stage cost function $\ell:\Re^n\times \Re^m\times \N\to\Re$
which satisfies Assumption~\ref{assum:well-posedness},
a pair $(x_s^\theta, u_s^\theta)$ is called an \textit{optimal steady state}
of~\eqref{eq:system_evo} subject to~\eqref{eq:constraints} 
with respect to $\ell$ if it is a minimizer of the problem
\begin{align*}
\ell_s(\theta) \dfn \min_{x,u} \left\{\ell(x,u, \theta){\mid}  f(x,u,\theta)=x, (x,u)\in Y_\theta\right\}
\end{align*}
\end{definition}
For reasons that will be better elucidated in the next section, we need to 
draw the following weak controllability assumption essentially requiring 
that if $x_k=x_s^i$ and $\theta_k=j$ then there is a control action $\bar{u}_{s}^{i,j}$ so 
that at time $k+1$ the state is steered to $x_{k+1}=x_s^{\bet(j)}$.
%
%
% ~~~~ ASSUMPTION
%      Weak controllability condition ~~~~
%
\begin{assumption}[Controllability]\label{assum:xs-invariant}
 In addition to Assumption~\ref{assum:well-posedness}, for all $i,j\in\NN$ there is a 
 control law $\bar{u}_s:\Re^n\times\NN\to\Re^m$ with $\bar{u}_s(x_s^i, j) = \bar{u}_s^{i,j}$
 so that $(x_s^i, \bar{u}_s^{i,j})\in Y_j$ and
 $f(x_s^i, \bar{u}_s^{i,j}, j) = x_s^{\bet(j)}$.
\end{assumption}

\subsection{Model predictive control}\label{sec:mpc-simple}

In this section we shall present a model predictive control 
framework for constrained Markovian switching systems with
mode-dependent optimal steady state points. 

Let $u_k\lhd \FF_k$ and $\uu_N=(u_0, \ldots, u_{N-1})$, 
and define $V_N$ 
 \begin{align*}
  V_N(x_0, \theta_0,\mathbf{u}_N) 
    &= \expect \Bigg[ V_f(x_{N}, \theta_{N})+\sum_{j=0}^{N-1}\ell(x_{j}, u_{j}, \theta_{j}) 
   \Big| \FF_0 \Bigg].
 \end{align*}
Here, we take $V_f=0$ and the state sequence 
satisfies~\eqref{eq:system_evo}.

We introduce the following stochastic economic model predictive control problem
%
%
% @@@@ PROBLEM #1 : Equality terminal constraints @@@@
%
\begin{subequations}\label{eq:mpc-00}
 \begin{align}
  \mathbb{P}(x,\theta)&: V_N^\star(x,\theta) = \inf_{\uu_N} V_{N}(x,\theta, \uu_N), \label{eq:MPC_sub1}
  \intertext{and for $k=0,\ldots, N-1$, subject to}
  x_{k+1} &= f(x_k, u_k, \theta_k) \label{eq:MPC_sub2}\\  
  (x_{k}, u_k) &\in Y_{\theta_k}\\
  (x_0, \theta_0)&= (x,\theta) \label{eq:MPC_sub3}\\
  x_N &= x_s^{\operatorname{bet}(\theta_{N-1})}\label{eq:bet}\\
  u_k &\lhd \FF_k.\label{eq:causality}
  \end{align}
  \label{eq:MPC_basic}
\end{subequations}
Because of Assumption~\ref{assum:well-posedness} and
in light of~\cite[Thm.~1.17]{rockafellar-book} the infimum in~\eqref{eq:mpc-00}
is attainable and the corresponding set of minimizers is compact.
Note that in the above formulation the minimization is carried out
in a space of control policies $\uu=\{u_0,\ldots, u_{N-1}\}$ where
$u_{k}$ are \textit{causal} control laws --- as required by~\eqref{eq:causality}.

Let $\uu^\star(x,\theta) = \{u_0^\star(x,\theta),  \ldots,
u_{N-1}^\star(x_{N-1},\theta_{N-1})\}$
be an optimizer of~\eqref{eq:MPC_basic}.
The receding horizon control law that accrues from this problem is
$\kappa_N(x,\theta)\dfn u_0^\star(x,\theta)$ and the closed-loop system
satisfies
\begin{align}\label{eq:mpc-controlled-system}
  x_{k+1} = f(x_k, \kappa_N(x_k, \theta_k), \theta_k).
 \end{align}

\subsection{Recursive feasibility}

We will now prove that the MPC problem in~\eqref{eq:MPC_basic}
is recursively feasible. 

%
%
% ------ PROPOSITION
%        Recursive feasibility of EMPC #1   ------
%
\begin{proposition}\label{prop:recursive-feasibility}
 Let $X_N\subseteq \Re^n \times \NN$ be the domain of problem $\mathbb{P}$.
 If Assumption~\ref{assum:xs-invariant} holds and problem $\mathbb{P}(x_k,\theta_k)$
 is feasible, then problem $\mathbb{P}(x_{k+1}, \theta_{k+1})$, with 
 $x_{k+1} = f(x_k, \kappa_N(x_k, \theta_k), \theta_k)$ and $\theta_{k+1}\in\mathcal{C}(\theta_k)$, is also feasible. 
\end{proposition}

\begin{pf}
 For given $(x,\theta)\in X_N$ let $\pi(x,\theta)=\{u_0^\star,\ldots,$
 $u_{N-1}\}$ be an optimizer of
 $\mathbb{P}(x,\theta)$ and let $x^\star(x,\theta)=\{x,
 x_1^{\star},\ldots, x_N^\star\}$ be the corresponding sequence 
 of states. Because of~\eqref{eq:bet} we have
 \[
  x_N^\star = x_s^{\bet(\theta_{N-1})}.
 \]
 Now take $x^+ = f(x,u_0^\star(x,\theta),\theta)$ and $\theta^+\in\mathcal{C}(\theta)$.
 We need to show that~$\mathbb{P}(x^+,\theta^+)$ is feasible. Take
 $\tilde{\pi}^+(x^+,\theta^+) \dfn \{u_1^\star,\ldots, u_{N-1}, u\}$ and let $u = \bar{u}_s(x_{N},\theta_{N})$.
 Then, by virtue of Assumption~\ref{assum:xs-invariant},
 $x^\star_{N+1} = x_s^{\bet(\theta_{N})}$, so $\tilde{\pi}^+$ will satisfy the 
 constraints of~$\mathbb{P}(x^+,\theta^+)$.~\hfill{$\Box$}
\end{pf}

\subsection{Performance assessment}

We will now prove that the closed-loop system has a bounded expected
asymptotic average cost (Theorem~\ref{thm:asymptotic-perf}).
First, we need to give the following result:

%
%
% ------ LEMMA
%        Necessary result for the performance bound   ------
%
\begin{lemma}\label{lem:lyap-type-ineq}
 Let Assumption~\ref{assum:xs-invariant} hold and let
\begin{align*}
\ell_N(\theta_k) {\dfn}
  \expect {\left[\ell(x_s^{\bet(\theta_{N-1})}, \bar{u}_s^{\bet(\theta_{N-1}),\theta_{N}},
  \theta_{N})\mid \theta_0 = \theta \right]},
\end{align*}
and  $\mathcal{L}V_N^\star(x_k, \theta_k) {\dfn} \expect[V_N^\star(x_{k+1}, \theta_{k+1}) -
      V_N^\star(x_k,\theta_k) {\mid} \FF_k]$; 
then, the following holds for all $(x_k, \theta_k)\in X_N$
 \begin{align}\label{eq:lyap-type-ineq}
  \mathcal{L}V_N^\star(x_k, \theta_k) \leq \ell_N(\theta_k) - \ell(x_k, \kappa_N(x_k,\theta_k), \theta_k).
 \end{align}
\end{lemma}

\begin{pf} % ----- Proof `lem:lyap-type-ineq` -----
 Let $(x,\theta)\in X_N$; then $\tilde\pi^+(x_{k+1}, \theta_{k+1})$ is feasible --- but not necessarily optimal ---
 for $\mathbb{P}(x_{k+1}, \theta_{k+1})$, therefore,
 $V_N^\star(x_{k+1},\theta_{k+1}) \leq V_{N}(x_{k+1},\theta_{k+1},\tilde{\pi}^+(x_{k+1}, \theta_{k+1}))$.
 By the tower property of the conditional expectation we know that
 $\expect[\expect[\cdot \mid \FF_{k+1}]\mid \FF_{k}] = \expect[\cdot \mid \FF_{k}]$
 since $\FF_k\subseteq \FF_{k+1}$.
 We then have
 \begin{align*}
  \mathcal{L}V_N^\star(x_k, \theta_k)
  &\leq \expect\Bigg[ \sum_{j=k+1}^{k+N-1}\ell(x_j,u_{j-k}^\star, \theta_j) + \ell(x_{k+N}, \bar{u}_s, \theta) - \\
  &\quad-\sum_{j=k}^{k+N-1}\ell(x_j, u_j, \theta_j)\mid \FF_k\Bigg]\\
  &=\expect \bigg[ \ell(x_s^{\bet(\theta_{k+N-1})}, \bar{u}_s^{\bet(\theta_{k+N-1}),\theta_{k+N}},\theta_{k+N})\\
  &\quad- \ell(x_k,\kappa_N(x_k,\theta_k),\theta_k) \mid \FF_k \bigg]\\
  &=\ell_N(\theta_k) -   \ell(x_k,\kappa_N(x_k,\theta_k),\theta_k),
 \end{align*}
 where $u_{j-k}^\star=u_{j-k}^\star(x_j, \theta_j)$ and this completes the proof.~\hfill{$\Box$}
\end{pf} % ----- Proof `lem:lyap-type-ineq` -----

The irreducibility and aperiodicity assumptions 
(Assumption~\ref{assum:well-posedness}) imply the existence 
of a limiting probability vector $\pi=(\pi^1,\ldots, 
\pi^\nu)\in\Re^\nu$ which satisfies $\pi P = \pi$ and does 
not depend on the initial distribution 
$v$~\citep{markov-chains-book}.
%
%
% ---- THEOREM: 
%      Asymptotic Performance Bound ----
%
\begin{theorem}[Asymptotic performance]\label{thm:asymptotic-perf}
 Let Assumption~\ref{assum:xs-invariant} hold and let $\{x_k\}_k$
 be a sequence satisfying~\eqref{eq:mpc-controlled-system}.
 Define the \textit{asymptotic average cost} as the random variable
 \begin{subequations}
 \begin{align}\label{eq:J}
  J \dfn \expect \bigg[ \limsup_{T\to \infty}\tfrac{1}{T}\sum_{k=0}^{T-1}\ell(x_k, u_k, \theta_k)\bigg]
 \end{align}
Then,
 \begin{align}\label{eq:as-perf-bounds}
  J \leq \ell_\infty \dfn \sum_{i \in \NN} \pi_i \ell_N (i).
 \end{align}
 \end{subequations}
\end{theorem}

\begin{pf} % ---- Proof of `thm:asymptotic-perf` ----
By taking asymptotic averages and the expectation with respect to $\FF_0$ 
on both sides of~\eqref{eq:lyap-type-ineq} we have
\begin{align}
  &\expect \left[\liminf_{T\to \infty}   \tfrac{1}{T}\sum_{k=0}^{T-1} \mathcal{L}V_N^\star(x_k, \theta_k) \right] \notag\\
  &\leq  \expect \bigg[ \liminf_{T\to \infty}   \tfrac{1}{T} \sum_{k=0}^{T-1} \ell_N(\theta_k) - \ell(x_k, \kappa_N(x_k,\theta_k), \theta_k)  \bigg]\notag\\
  &\leq\expect \bigg[\liminf_{T\to \infty}  \tfrac{1}{T} \sum_{k=0}^{T-1}\ell_N(\theta_k) \notag\\
  & \quad-  \limsup_{T\to \infty} \tfrac{1}{T}\sum_{k=0}^{T-1} \ell(x_k,\kappa_N(x_k,\theta_k),\theta_k)  \bigg]\notag\\
  &\leq\liminf_{T\to \infty}  \expect \bigg[\tfrac{1}{T} \sum_{k=0}^{T-1}\ell_N(\theta_k) \bigg] \notag\\
  & \quad-  \expect \bigg[ \limsup_{T\to \infty} \tfrac{1}{T}\sum_{k=0}^{T-1} \ell(x_k,\kappa_N(x_k,\theta_k),\theta_k)  \bigg].\label{eq:proof-bound-x}
  \end{align}
We now use the fact that $\expect [\ell_N(\theta_k)] 
= \sum_{i\in\NN}\pi_k^i \ell_N(i)$, where $\pi_k^i = 
\prob[\theta_k = i]$ and since $\pi_k^i \to \pi^i$
as $k\to\infty$, 
we have that $\expect[\ell_N(\theta_k)]\to \ell_\infty$
and the right hand side of~\eqref{eq:proof-bound-x} is equal 
to $\ell_\infty - J$.

Using \cite[Lemma~19]{PatSopSarBem14} and because of the fact 
that $\ell$ are nonnegative,
\begin{align*}
  &\expect \left[
	\liminf_{T\to \infty}  
	  \tfrac{1}{T}\sum_{k=0}^{T-1} 
	    \mathcal{L}V_N^\star(x_k, \theta_k) 
      \right] \\
  & = \expect \left[
	  \liminf_{T\to \infty}
	    \tfrac{1}{T} (V_N^\star(x_{T}, \theta_{T}) 
	  - V_N^\star(x_0,\theta_0))  
	 \right]\\
  & \ge \liminf_{T\to \infty} 
	\left( 
	  -\tfrac{1}{T} V_N^\star(x_0,\theta_0) 
	\right) = 0.
 \end{align*}
Combining the two results gives
\begin{align*}
\expect \left[\limsup_{T\to \infty}\tfrac{1}{T}\sum_{k=0}^{T-1} \ell(x_k,\kappa_N(x_k,\theta_k),\theta_k)  \right] \le \ell_\infty
\end{align*}
which completes the proof.~\hfill{$\Box$}
\end{pf}  % ---- Proof of `thm:asymptotic-perf` ----

\subsection{Mean square stability}
We will now study under what conditions a
Markovian system is mean square stable towards an equilibrium point.
%
%
% ---- ASSUMPTION
%      Common equilibrium ----
% 
\begin{assumption}[Common optimal equilibrium]\label{assum:common_point}
There exists one common optimal stationary point $(x_s,u_s)$ for all modes
which is the solution of the optimization problem in Definition~\ref{def:optimal-steady-states}
and, without loss of generality, assume $x_s=0$, $u_s=0$.
\end{assumption}

Consider the following Markovian switching system
\begin{align}\label{eq:closed-loop-markovian}
 x_{k+1} = f(x_k, \theta_k),
\end{align}
and let $r_k=(\theta_0,\ldots,\theta_k)$ be an admissible \textit{switching sequence} starting
from $\theta_0$ and $\phi(k;x_0,r_k)$ be the trajectory of~\eqref{eq:closed-loop-markovian}
with $\phi(0;x_0, r_0)=x_0$. We recall the definition of mean square stability
\begin{definition}[Mean Square Stability]
 We say that~\eqref{eq:closed-loop-markovian} is \textit{mean square 
 stable} if $\expect[\|\phi(k;x_0, r_k)\|^2] \to 0$,
as $k\to \infty$ for all $x_0$ and $\theta_0$.
\end{definition}
We extend the notion of dissipativity to Markovian systems 
as follows
%
%
% ---- DEFINITION: Stochastic dissipativity ----
%
\begin{definition}[Stochastic dissipativity]
\label{def:stochastic-dissipativity}
% STOCHASTIC DISSIPATIVITY
We say that system~\eqref{eq:closed-loop-markovian} is 
\textit{stochastically dissipative} with respect to a 
stochastic supply rate $s:\Re^n\times \Re^m\times \NN\to \Re$ 
if there is a function $\lambda:\Re^n\times \NN
\to\Re$, lower semicontinuous in the first argument, so that for all $x_k\in\Re^n$ and $\theta_k \in \NN$
\begin{align}
  \mathcal{L}\lambda(x_k,\theta_k)\leq s(x_k, u_k, \theta_k).
 \label{eq:dissipativity_basic}
\end{align}
where $\mathcal{L}\lambda(x_k,\theta_k) \dfn \expect [\lambda(x_{k+1}, \theta_{k+1}) - \lambda(x_k, \theta_k) \mid \mathfrak{F}_k ]$.
% STRICT STOCHASTIC DISSIPATIVITY
We say that~\eqref{eq:system_evo} is \textit{strictly 
stochastically dissipative} with respect to $s$ if there is a 
convex function $\rho:\Re^n\times \NN \to \Re_+$, 
positive definite with respect to $x_s$,  so that the left 
hand side of~\eqref{eq:dissipativity_basic} is no larger than
$s(x_k, u_k, \theta_k) - \rho(x_k,\theta_k)$.
\end{definition}

%
% ..... ASSUMPTION: Strict stoch. dissipativity .....
%
\begin{assumption}[Strict stochastic dissipativity]
\label{assum:strict_dissipativity}
  Function $\lambda(x_s, \theta)$ is independent of $\theta$
  and let $\lambda_s \dfn \lambda(x_s, \theta)$.
  In addition to Assumption~\ref{assum:common_point}, 
  system~\eqref{eq:closed-loop-markovian} is 
  strictly stochastically dissipative with storage function 
  $s(x,u,\theta) = \ell(x,u,\theta) -\ell_s$.
\end{assumption}

%
% ---- DEFINITION: Rotated stage cost ----
%
Let us define the \emph{rotated stage cost} function as
\begin{align}
  L(x_{k},u_{k},\theta_{k}) &\dfn 
     \ell(x_{k},u_{k},\theta_{k}) 
   - \mathcal{L}\lambda(x_k,\theta_k).
  \label{eq:rotated_steg_cost}
\end{align}
We now define the \textit{rotated cost function} 
$\widetilde{V}_N(x,\theta,\mathbf{u}_N)$ as follows
\begin{align*}
  \widetilde{V}_N(x_0, \theta_0,\mathbf{u}_N) 
  &= \expect \Bigg[ 
     \sum_{j=0}^{N-1}L(x_{j}, u_{j}, \theta_{j})  
     \Big| \FF_0 
     \Bigg]
 \end{align*}
using again $V_f = 0$ 
and we introduce the rotated MPC problem
 \begin{align}\label{eq:mpc-rotated}
  \bar{\mathbb{P}}(x,\theta)&: 
  \widetilde{V}_N^\star(x,\theta) 
    = \inf_{\uu_N} \widetilde{V}_{N}(x,\theta, \uu_N),
  \end{align}
subject to~\eqref{eq:MPC_sub2}--\eqref{eq:causality}.

%
% ---        LEMMA: LV_N^* <= -rho       ----
% +++ using this lemma we prove stability +++
%
\begin{lemma}\label{lemma:lyapunov-ineq}
 Problem~$\bar{\mathbb{P}}(x,\theta)$ is recursively feasible 
 and it has the same set of minimizers as~$\mathbb{P}(x,
 \theta)$.
 Let $\tilde{\kappa}_N$ be the receding horizon control law 
 which accrues from $\bar{\mathbb{P}}(x,\theta)$.
 If Assumption~\ref{assum:strict_dissipativity} holds, then
 \begin{align}
   \mathcal{L}\widetilde{V}_N^\star(x_k, \theta_k) 
  \le - \rho(x_{k},\theta_{k}),
 \end{align}
 where $\rho\Re^n\times \NN \to \Re_+$ is
 a positive definite function in the first argument 
 with respect to $x_s$.
\end{lemma}

\begin{pf} % ----- lemma `lemma:lyapunov-ineq` -----
Problems $\mathbb{P}$ and $\bar{\mathbb{P}}$ have the same set of constraints,
therefore, they have the same feasibility domain and the recursive feasibility of
$\bar{\mathbb{P}}$ follows from Proposition~\ref{prop:recursive-feasibility}.
Rotated cost function can be expanded as
\begin{align*}
  &\widetilde{V}_N(x_k,\theta_k,\mathbf{u}_N) = \expect \Bigg[\sum_{j=k}^{k+N-1} L(x_j,u_j, \theta_j) \Big| \mathfrak{F}_k \Bigg]   \\
  &= \expect \Bigg[
    \sum_{j=k}^{k+N-1} 
       \ell(x_j,u_j, \theta_j)  
     - \mathcal{L}\lambda(x_k, \theta_k)
    \Big| 
    \mathfrak{F}_k 
    \Bigg] 
\end{align*}
We now use the fact that
\begin{align*} 
&\expect\left[\sum_{j=k}^{k+N-1} \mathcal{L}\lambda(x_k, \theta_k)\mid \FF_{k}\right] \\
  = &\expect[\lambda(x_{k+N-1},\theta_{k+N-1})
  -\lambda(x_k, \theta_k)\mid \FF_k] \\
  = &\lambda_s - \lambda(x_k, \theta_k).
\end{align*}
Therefore,
\begin{align*}
 \widetilde{V}_N(x_k,\theta_k,\mathbf{u}_N) 
 = V_N(x_k,\theta_k,\mathbf{u}_N)  + \lambda(x_{k}, \theta_{k})-\lambda_s.
\end{align*}

The rotated and original cost functions differ only by a 
constant so the two problems, $\mathbb{P}$ and 
$\bar{\mathbb{P}}$, share a common optimal sequence.
Proceeding as in Lemma~\ref{lem:lyap-type-ineq} the following holds
\begin{align}
 \mathcal{L}\tilde{V}_N^\star(x_k,\theta_k)  
 \le \ell_s -L(x_k,\tilde{\kappa}_N(x_k,\theta_k),\theta_k),
\end{align}
By tracing the arguments of~\cite{rawlings2012fundamentals},
$L(x_k,u_k,\cdot) \ge \ell_s$.
Combining~\eqref{eq:rotated_steg_cost} and Assumption~\ref{assum:strict_dissipativity}  we arrive at
\begin{align}
 L(x_k,u_k,\theta_k) \ge \rho(x_k,\theta_k) + \ell_s,
\end{align}
which completes the proof.~\hfill$\Box$ 
\end{pf} % ----- lemma `lemma:lyapunov-ineq` -----

Next, we draw an additional assumption on $\rho(\cdot,\theta)$:
%
% ---- ASSUMPTION: 
%      `assum:lower-quad-bound`
%      Quadratic lower bound on \rho ----
%
\begin{assumption}[Quadratic lower bound]
\label{assum:lower-quad-bound}
 There exist a positive constant $\gamma$, such that
$% $ \begin{align}
    \rho(x,i) \ge \gamma \| x - x_s \|^2$
 holds for all $x$.
\end{assumption}

%
%
% ~~~~ THEOREM: 
%      Main MS stability result ~~~~
%
\begin{theorem}
 Suppose Assumption~\ref{assum:lower-quad-bound} is satisfied. Then,
 system~\eqref{eq:closed-loop-markovian} is MSS.
\end{theorem}

\begin{pf}
 All assumptions required by~\citep[Theorem 24]{PatSopSarBem14} are met and 
 entail mean square stability.~\hfill{$\Box$}
\end{pf}

\section{Uniform Invariance and Terminal Constraints}
In this section we relax the restrictive requirement 
$x_N = x_s^{\bet(\theta_{N-1})}$
and we instead replace it with a terminal constraint 
of the form $(x_N, \theta_N)\in X^f$
along with a terminal penalty function $V_f$ and we 
derive conditions so that the controlled system is 
mean-square stable.

We will now make use of the following definition~\citep{PatSopSarBem14}
%
% ++++ DEFINITION: 
%      Uniform positive invariance ++++
%
\begin{definition}[Uniform positive invariance]
A family of nonempty sets $C = \{C_i\}_{i\in\NN}$ is said
to be \emph{uniformly positive invariant} (UPI) for the constrained Markovian
switching system~\eqref{eq:closed-loop-markovian} if for every $x_k\in C_{\theta_k}$, $x_{k+1}\in C_{\theta_{k+1}}$.
\end{definition}

As before, we assume that there is one stationary point $\ell_s$ and require, with a slight abuse of notation, that  $\lambda_s = \lambda(x_s,\theta), V_f(x_s) = V_f(x_s,\theta)$
for all $\theta \in \NN$.
Now we make a central assumption regarding our exposition

%
%
% ----- ASSUMPTION:
%       Terminal control region
%       Terminal controller      -----
%
\begin{assumption}[Terminal control law]
\label{assum:terminal-region-law}
There exists a control law $\kappa_f:\Re^n\times \NN\to \Re^m$ and a 
collection of sets $X^f=\{X^f_i\}_{i\in\NN}$ so that 
\begin{itemize}
 \item[i.] $X^f$ is UPI for the closed-loop system controlled by $\kappa_f$ and
 \item[ii.] for  all $(x,\theta)\in X^f$
\begin{align}\label{eq:terminal-region-law}
	\mathcal{L}V_f(x_k, \theta_l) \leq  -\ell(x_k,\kappa_f(x_k,\theta_k),\theta_k) + \ell_s.
\end{align}
\end{itemize}
\end{assumption}

Now consider the following stochastic economic model predictive
control problem

%
%
% @@@@ PROBLEM #2 : Terminal constraints @@@@
%
\begin{subequations}\label{eq:mpc-terminal}
 \begin{align}
  \mathbb{P}_T(x,\theta)&: V_N^\star(x,\theta) = \inf_{\uu_N} V_{N}(x,\theta, \uu_N)\label{eq:MPC_ocp_T}
  \end{align}
and for $k=0,\ldots, N-1$, it is subject to  
  \begin{align}
  x_{k+1} &= f(x_k, u_k, \theta_k) \label{eq:MPC_ocp_T_sub1}\\
  (x_{k}, u_k) &\in Y_{\theta_k} \label{eq:MPC_ocp_T_sub2}\\
  (x_0, \theta_0)&= (x,\theta) \label{eq:MPC_ocp_T_sub2b} \\
  x_N &\in X^f_{\theta_k} \label{eq:set-constraint}\\
  u_k &\lhd \mathfrak{F}_k. \label{eq:MPC_ocp_T_sub3}
  \end{align}
  \label{eq:MPC_terminal_set}
\end{subequations}

Again, the same reasoning as in section~\ref{sec:mpc-simple} 
applies regarding the existence of optimal solutions.
Let $\hat{\uu}^\star(x,\theta) = \{u_0^\star(x,\theta), 
\ldots, u_{N-1}^\star(x_{N-1},\theta_{N-1})\}$
be an optimizer of~\eqref{eq:mpc-terminal}.
The receding horizon control law is given by 
$\hat{\kappa}_N(x,\theta) \dfn u_0^\star(x,\theta)$. 

In light of the state-input constraints~\eqref{eq:MPC_ocp_T_sub2}
we must require that the sets $X^f_i$ in Assumption~\ref{assum:terminal-region-law}
are subsets of $X_N$, the feasibility domain of $\mathbb{P}_T$.

\subsection{Recursive feasibility}
Here, we will show that stochastic economic model predictive
control problem~\eqref{eq:mpc-terminal} is recursively feasible.

%
% [[ Proposition ]]
% The problem with terminal constraints is
% recursively feasible.
%
\begin{proposition}
\label{prop:recursive-feasibility-terminal-set}
 Let $X_N\subseteq\Re^n\times \NN$ be the feasibility domain of $\mathbb{P}_T$ and 
 let Assumption~\ref{assum:terminal-region-law}-i
 hold. Then, $X_N$ is UPI for the MPC-controlled system.
\end{proposition}

\begin{pf} % ----- Proof `prop:recursive-feasibility-terminal-set` -----
 For given $(x,\theta)\in X_N$ let $\pi(x,\theta)=\{u_0,\ldots,u_{N-1}^\star\}$ be an optimizer of
 $\mathbb{P}_T(x,\theta)$ and let $x^\star(x,\theta)=\{x, x_1^\star,\ldots,$ $x_N^\star\}$ be the corresponding sequence of
 states.
 Because of~\eqref{eq:set-constraint} we have
 $x_N^\star \in X^f_{\theta_N}$.
 Now take $x^+ = f(x,u_0^\star(x,\theta),\theta)$, $\theta^+\in\mathcal{C}(\theta)$ and let
 $\tilde{\pi}^+(x^+,\theta^+) \dfn \{u_1^\star,\ldots, u_{N-1}^\star, u_f\}$
 where $u_f= \kappa_f(x_{N},\theta_{N})$.
 Then, since $X^f$ is a UPI set,
 $(x_{N+1},\theta_{N+1}) \in X^f$, so $\tilde{\pi}^+$ satisfies the constraints 
 of~$\mathbb{P}_T(x^+,\theta^+)$.~\hfill{$\Box$}
\end{pf} % ----- Proof `prop:recursive-feasibility-terminal-set` -----

\subsection{Expected asymptotic average performance}
Here we show that the asymptotic average cost of the EMPC-controlled system with
terminal constraints is no higher than the cost of the  best stationary point.
%
%
% ---- THEOREM:
%      Asymptotic performance bound
%      J <= ell_s                     ----
%
\begin{theorem}
Let Assumption~\ref{assum:terminal-region-law} hold and let $\{x_k\}_k$
be a sequence satisfying~\eqref{eq:mpc-controlled-system} with $u_k = \hat{\kappa}_N(x_k, \theta_k)$.
 Then, $ J  \leq \ell_s$.
\end{theorem}

\begin{pf}
Using the optimal solution $\pi(x,\theta)$ of~\eqref{eq:mpc-terminal} with initial conditions $(x,\theta)$
we construct a feasible shifted policy $\pi^+(x^+,\theta^+)$
as in the proof of the Proposition~\ref{prop:recursive-feasibility-terminal-set}. 
Then $V_N^\star(x_{k+1},\theta_{k+1}) \leq V_N(x^{+}, \tilde{\pi}^+,\theta^+)$ and
\begin{align}
    &\mathcal{L}V_N^\star(x_k,\theta_k)= \expect\Bigg[ \sum_{j = k + 1}^{k+N-1}\ell(x_j,u_j^\star, \theta_j)   \nonumber\\ 
    & + \ell(x_{k+N}, \kappa_f(x_{k+N},\theta_{k+N}), \theta_{k+N}) + V_f(x_{k+N+1}, \theta_{k+N+1}) \nonumber \\
    & -\sum_{j=k}^{k+N-1}\ell(x_j, u_j^{\star}, \theta_j) - V_f(x_{k+N}, \theta_{k+N})\mid \mathfrak{F}_k\Bigg] \nonumber \\
    &  \leq \ell_s - \ell(x,\hat{\kappa}_N(x,\theta),\theta).  \nonumber
\end{align}
Here, we used tower property and Assumption~\ref{assum:terminal-region-law}. Proceeding as in 
Theorem~\ref{thm:asymptotic-perf} we prove the assertion.~\hfill{$\Box$}
\end{pf}

\subsection{Mean square stability}
In this section we will give conditions under which Markovian system with terminal region constraint is mean square stable 
towards a common equilibrium point. Once again, our main argument will be the equivalence between original and suitably \emph{rotated} problem.

We define the following \emph{rotated terminal function}
\begin{align}
	\tilde{V}_f (x_k,\theta_k) = V_f (x_k,\theta_k) + \lambda(x_{k}, \theta_{k}) - V_f (x_s) - \lambda_s.
\end{align}

Combining condition~\eqref{eq:dissipativity_basic} (Definition~\ref{def:stochastic-dissipativity}) with the rotated stage cost we may easily derive
\begin{align}
	L(x_k,u_k,\theta_k) \ge \rho(x_k, \theta_k).
	\label{eq:stage_cost_rho}
\end{align}
%

%
% ----- LEMMA
%       LVf <= -L (rotated) -----
%
\begin{lemma}\label{lemma:lyapunov-type-rotated}
 Suppose Assumption~\ref{assum:terminal-region-law} holds. Then
\begin{align}
	\mathcal{L}\tilde{V}_f(x_k, \theta_k) \leq  -L(x_k,\kappa_f(x_k, \theta_k),\theta_k).
	\label{eq:terminal_decrease_assumption_rotated}
\end{align} 
\end{lemma}

\begin{pf} % ----- Proof `lemma:lyapunov-type-rotated` -----
We add $\mathcal{L}\lambda(x_k,\theta_k)$ to both sides of~\eqref{eq:terminal-region-law}
\begin{subequations}
 \begin{align*}
  &\mathcal{L}\tilde{V}_f(x_k, \theta_k)+\mathcal{L}\lambda(x_k,\theta_k)
  \leq - \ell(x_k,\kappa_f(x_k, \theta_k),\theta_k) + \ell_s\\
  &\quad+ \mathcal{L}\lambda(x_k,\theta_k).
  \end{align*}
\end{subequations}
The right hand side is equal to the rotated stage cost
 \begin{align*}
  &\expect \bigg[ V_f(f(x_k,\kappa_f(x_k,\theta_k)),\theta_{k+1}) + \lambda(x_{k+1}, \theta_{k+1}) \\
  &- V_f( x_k, \theta_{k}) - \lambda(x_k, \theta_k) \mid \FF_{k} \bigg] \leq -L(x_k,\kappa_f(x_k, \theta_k),\theta_k) \nonumber.
  \end{align*}
We add  $V_f (x_s) + \lambda_s - V_f (x_s) - \lambda_{s}$ 
to the left hand side and, after rearranging, arrive 
at~\eqref{eq:terminal_decrease_assumption_rotated}.~\hfill{$\Box$}
\end{pf} % ----- Proof `lemma:lyapunov-type-rotated` -----

Now, we introduce a \emph{rotated stochastic economic MPC problem}
 \begin{align}\label{eq:mpc-rotated-terminal}
  \bar{\mathbb{P}}_T(x,\theta)&: \tilde{V}_N^\star(x,\theta) 
  = \inf_{\uu_N} \tilde{V}_{N}(x,\theta, \uu_N)
  \end{align}
subject to~\eqref{eq:MPC_ocp_T_sub1}-\eqref{eq:MPC_ocp_T_sub3}.

% 
% ----- THEOREM 
%       Equivalence of the two problems  -----
%
\begin{theorem}
 Problem~$\bar{\mathbb{P}}_T(x,\theta)$ is recursively feasible 
 and has the same set of minimizers as~$\mathbb{P}_T(x,\theta)$.
\end{theorem}

\begin{pf} % ----- Proof: equivalence of problems -----
Problems $\mathbb{P}_T$ and $\bar{\mathbb{P}}_T$ have the same set 
of constraints, therefore, they have the same feasibility domains 
and the recursive feasibility of $\bar{\mathbb{P}}$ follows from 
Proposition~\ref{prop:recursive-feasibility-terminal-set}.
The rotated cost function can be expanded as
$
  \tilde{V}_N(x_k,\theta_k,\mathbf{u}_k) = \expect [\sum_{j=k}^{k+N-1} L(x_j,u_j, \theta_j) + \tilde{V}_f(x_N,u_N,\theta_N) \mid \mathfrak{F}_k ]  
  = \expect [\sum_{j=k}^{N-1} (\ell(x_j,u_j, \theta_j) + \lambda(x_j,\theta_j)
  - \expect [\lambda(x_{j+1},\theta_{j+1} - \ell_s)\mid \mathfrak{F}_j]) + V_f (x_N,\theta_N) + \lambda(x_{N}, \theta_{N})  - V_f (x_s) - \lambda_s) \mid \mathfrak{F}_N ] \mid \mathfrak{F}_k]      
  = V_N(x,\mathbf{u},\theta)  + \lambda(x,\theta) - N \ell_s - V_f(x_s) -  \lambda_s.
$
The two cost functions, $V_N$ and $\tilde{V}_N$ differ by 
feedback-invariant quantities, hence, the optimal solutions
of the two problems will coincide.~\hfill{$\Box$}
\end{pf} % ----- Proof: equivalence of problems -----

\begin{theorem}\label{thm:mss_of_terminal_set_MPC} 
 Suppose Assumptions~\ref{assum:lower-quad-bound} and~\ref{assum:terminal-region-law} are satisfied. Then,
 system~\eqref{eq:mpc-controlled-system} is MSS
 with domain of attraction $X_N$.
\end{theorem}
\begin{pf}
 All assumptions required by~\citep[Theorem 24]{PatSopSarBem14} are met and 
 we can infer mean square stability.~\hfill{$\Box$}
\end{pf}

\subsection{Linearization-based design}
In this section we demonstrate how to design a terminal cost function and give a terminal control 
law using local linearization around origin.  In other words, we give conditions under which 
Assumption~\ref{assum:terminal-region-law}\,-ii is satisfied, given that  
Assumption~\ref{assum:terminal-region-law}\,-i holds for a nonlinear system with a particular 
control law. In the next section we shall also demonstrate how to design an ellipsoidal set 
$X^f$ such that it satisfies Assumption~\ref{assum:terminal-region-law}\,-i.

To simplify the notation let
$\bar{\ell}(x,\theta) = \ell(x,\kappa_f(x,\theta),\theta) - \ell(0,0,\theta) $ for all $\theta \in \NN$, 
be a shifted stage cost function. Define $\hat{f}_\theta(x) \dfn f(x,\kappa_f(x,\theta),\theta)$, 
where $\kappa_f(x,\theta)$ is a terminal control law that we will introduce shortly.
The evolution of the nonlinear system is described by $x_{k+1} = \hat{f}_\theta(x_k)$, for all $\theta \in \NN$.

To proceed we need the following assumption which is weaker 
than twice differentiability which is commonly used in 
the literature~\citep{RawlingsMayne-book}.
%
%
% <---- ASSUMPTION 
%       `assum:beta_smooth_system`
%       beta-smoothness             ----->
%
\begin{assumption}[Smoothness]\label{assum:beta_smooth_system}
  Functions $\hat{f}_\theta(x)$ are $\beta_f^\theta$-smooth and 
  $\bar{\ell}(x,\theta)$ are $\beta_\ell^\theta$-smooth for all $\theta \in \NN$.
\end{assumption}

Let 
\begin{align}
  z_{k+1} = A_{\theta_k} z_k + B_{\theta_k} u_k
  \label{eq:linearized_MJLS}
\end{align}
be the corresponding linearized Markovian jump linear systems (MJLS), where
$A_i = \frac{\partial f_i}{\partial x} (0,0)$ and $B_i = \frac{\partial f_i}{\partial u} (0,0)$ for all $i \in \NN.$
Hereafter, we will make the following assumption:

%
%
% (---- ASSUMPTION 
%       MS-stabilizability -----)
%
\begin{assumption}
The set of pairs $\{(A_i,B_i)\}_{i \in \NN}$ is mean square stabilizable.
\label{ass:stabilizable_linearization}
\end{assumption}

\cite{mjls2005costa} provide conditions for Assumption~\ref{ass:stabilizable_linearization}
to hold. We recall the following result for MJLS~\citep{PatSopSarBem14}

%
%
% ~~~~~ PROPOSITION
%       `prop:mss-for-mjls`
%       MSS for MJLS         ~~~~~
%
\begin{proposition}[MSS of MJLS]
\label{prop:mss-for-mjls}
Consider system~\eqref{eq:linearized_MJLS} subject to~\eqref{eq:constraints}
in closed loop with $\kappa(x,i) = K_i x$. Suppose there is a UPI set $X^f$
and matrices $P^f = \{P^f_i\}_{i\in\NN}$ so that $P^f_i \succcurlyeq \Gamma_i^\top \EE_i(P^f)\Gamma_i + Q^*_i$
with $\Gamma_i \dfn A_i + B_i K_i$, $\EE_i(P^f) \dfn \sum_{j \in \mathcal{C}(i)} p_{ij}P_j^f$
and $Q^*_i=(Q^*_i)^\top\succ 0$. Then, the closed-loop system is MS stable in $X^f$.
\end{proposition}

Next, we will design a terminal cost function 
$V_f(x,\theta)$ which, under certain assumptions 
(see Theorem~\ref{thm:lin-stability}) 
satisfies a desired Lyapunov-type inequality (see 
Assumption~\ref{assum:terminal-region-law}\,-ii).

First, we design a quadratic cost function 
$\ell_q(x,\theta)$ which is an upper bound on the shifted cost.
%
%
% ~~~~~ LEMMA
%       `lemma:l_q`
%       Useful result to prove MSS for the 
%       linearized system                   ~~~~~
%
\begin{lemma}
 Let $\ell_q(x,\theta) \dfn \tfrac{1}{2}x^\top Q^*_{\theta}x + q_{\theta}^\top x$ where $Q^*_\theta = (\alpha+\beta_\ell^\theta) I,
 q_{\theta} = \nabla \bar{\ell}(0,\theta)$. Then it holds that $\ell_q(x,\theta) \ge \bar{\ell}(x,\theta) + \frac{\alpha}{2}\|x\|^2 $ for any $\alpha > 0$, for all $\theta\in\NN$.
 \label{lemma:l_q}
\end{lemma}
\begin{pf} % ------ Proof `lemma:l_q` -------
By Assumption~\ref{assum:beta_smooth_system} on 
$\bar{\ell}(x,\theta)$, we have that
$| \bar{\ell}(x,\theta) - q_{\theta}^\top x | 
\le  \beta^{\theta}_l / 2 \| x \|^2$. 
Adding $\alpha / 2 \| x \|^2$ to both sides the assertion follows.~\hfill{$\Box$}
\end{pf} % ------ Proof `lemma:l_q` -------

We may now choose our terminal cost to be the following infinite sum
\begin{align}
 V_f(x,i) = \expect \left [ \sum_{k=0}^{\infty} \ell_q\left(x_k,\theta_k \right)\Big|\ \FF_0  \right ],
 \label{eq:V_f}
\end{align}
for the MJLS  $x_{k+1} = \Gamma_{\theta_k} x_k$, with $x_0 = x, \theta_0 = \theta.$

Using the linearity of expectation we have
$V_f(x,\theta) = \expect \left [ 
\sum_{k = 0}^{\infty} 
  \tfrac{1}{2} x_k^\top Q_{\theta_k}^{*}x_k 
  \right]
 + \expect \left[ 
   \sum_{k = 0}^{\infty} 
     q_{\theta_k}^\top x_k 
\right ]$ and $V_f$ can be written 
in the form
\begin{align}
 V_f(x,i) =  \tfrac{1}{2}x^\top P_i^{f} x + p_i^\top x,
\end{align}
where $P^f_i$ are computed as in Prop.~\ref{prop:mss-for-mjls} with $=$ 
in lieu of $\succcurlyeq$~\citep[Prop.~3.20]{mjls2005costa}.
Because of the parametrization of  $Q^*_i$ in Lemma~\ref{lemma:l_q}, we may choose 
$P_i^f = P_i^{\beta} + \alpha P_i^I$
and require that 
\begin{subequations}
 \begin{align}
         P^I_i  &= I + \Gamma_i^\top \EE_i( P^I ) \Gamma_i,\\
     P_i^\beta  &= \beta_\ell^i I + \Gamma_i^\top \EE_i( P^{\beta} ) \Gamma_i 
 \end{align}
\end{subequations}

For convenience we re-introduce operator $\mathcal{L}$, 
but this time with a distinction between nonlinear and 
linear systems:
\begin{itemize}
 \item[i. ] $\mathcal{L}V_f(x_k,\theta_k) 
 = \expect [ 
    V_f ( \hat{f}_{\theta_k}(x_k),\theta_{k+1} ) 
  - V_f(x_k,\theta_k) \mid \mathfrak{F}_k]$ 
 \item[ii. ] $\mathcal{L}V_f^{\mathrm{lin}}(x_k,\theta_k) 
           = \expect [ 
             V_f \left( \Gamma_{\theta_k} x_k,\theta_{k+1} \right) 
           - V_f(x_k,\theta_k) \mid \mathfrak{F}_k]$.
\end{itemize}

Parameter $\alpha$ will be used to bound the mismatch between
$\mathcal{L}V_f(x_k,\theta_k)$ and $\mathcal{L}V_f^{\mathrm{lin}}(x_k,\theta_k)$ 
and a method for choosing it is presented in the proof of the next theorem.

%
%
% ....... THEOREM
%         Linearization-based MS stability
%         thm:lin-stability                 ........
%
% 
\begin{theorem}\label{thm:lin-stability}
 Consider the control law $\kappa_f(x,i) = K_ix$ and let 
 Assumptions~\ref{assum:beta_smooth_system} 
 and~\ref{ass:stabilizable_linearization} hold.
 Then $\mathcal{L}V_f(x,\theta) \leq -\bar{\ell}(x, \theta)$
 for $x\in\mathcal{B}_\delta$ for some $\delta > 0$. If 
 $X^f$ satisfies Assumption~\ref{assum:terminal-region-law}\,-i 
 with $X^f_i \subseteq \mathcal{B}_\delta$ and 
 Assumption~\ref{assum:lower-quad-bound} is satisfied,
 the controlled system is locally mean square stable.
\end{theorem}

\begin{pf} % ------ Proof `thm:lin-stability` ------
Let us introduce the  shorthand 
$ 
  \Delta \mathcal{L} V_f(x_k,\theta_k) 
    \dfn \expect [  
         V_f ( \hat{f}_{\theta_k}(x_k),\theta_{k+1} ) 
       - V_f \left( \Gamma_{\theta_k} x_k,\theta_{k+1} \right)  
      \mid \FF_k].
$
By the linearity of the conditional expectation we have 
$\mathcal{L}V_f(x,\theta) = \mathcal{L}V_f^{\text{lin}}(x,\theta) + \Delta \mathcal{L} V_f(x,\theta)$.
Because of~\eqref{eq:V_f}, the first term is $\mathcal{L}V_f^{\text{lin}}(x,\theta)  = -\ell_q(x,\theta)$.
The last term, after introducing $e(x,\theta) \dfn \hat{f}_{\theta}(x) - \Gamma_{\theta} x$, amounts to
\begin{align}
	&\Delta \mathcal{L} V_f(x,\theta)  = \tfrac{1}{2} e(x,\theta)^\top \EE_\theta(P^f)e(x,\theta)  \nonumber \\
	&\ - (\Gamma_{\theta}x)^\top \EE_\theta(P^f) e(x,\theta) + \EE_\theta(p)^\top  e(x,\theta).
	\label{eq:linearization_error_e}
\end{align}
where $e(x,\theta)$ is the linearization error. Under Assumption~\ref{assum:beta_smooth_system}
$\left\Vert e(x,\theta) \right\Vert \leq \tfrac{\beta^\theta_f}{2} \left\Vert x  \right\Vert^2$, therefore,
{
\begin{align}
	\Delta &\mathcal{L} V_f(x,\theta) \le \tfrac{(\beta^{\theta}_f)^2}{8} \| \EE_\theta(P^f) \| \| x \|^4 \nonumber \\
	&+\tfrac{\beta^{\theta}_f}{2} \| \Gamma_{\theta}\| \|\EE_\theta(P^f)\|\|x\|^3  + \tfrac{\beta^{\theta}_f}{2}\|\EE_\theta(p) \| \|x\|^2.
\end{align}
We need to show that $\Delta \mathcal{L} V_f(x,\theta)$ is upper bounded
by $\tfrac{\alpha}{2}\|x\|^2$ in a region of the origin for adequately large $\alpha$.
Recall that $\EE_\theta(P^f)$ depends on $\alpha$ as follows
\begin{align}
 \EE_\theta(P^f) = \EE_\theta(P^\beta) + \alpha \EE_\theta(P^I).
\end{align}
Using the triangle inequality
\begin{align}
  &\Delta\mathcal{L} V_f(x,\theta) \le \tfrac{(\beta^{\theta}_f)^2}{8} \|\EE_\theta(P^\beta)\| \|x\|^4 \nonumber \\
	&+\tfrac{\beta^{\theta}_f}{2} \| \Gamma_{\theta}\| \|\EE_\theta(P^\beta)\| \|x\|^3  + \tfrac{\beta^{\theta}_f}{2}\|\EE_\theta(p) \| \| x \|^2\notag\\
	&+\alpha\left(
	\tfrac{(\beta^{\theta}_f)^2}{8} \| \EE_\theta(P^I) \| \| x \|^4 
	+ \tfrac{\beta^{\theta}_f}{2} \| \Gamma_{\theta} \| \| \EE_\theta(P^I)  \| \| x\|^3 
	\right)
\end{align}
For the right hand side of the last inequality to be upper bounded by
$\tfrac{\alpha}{2}\|x\|^2$ it suffices to take $x\in\mathcal{B}_\delta$
with $\delta > 0$ and 
\begin{align*}
 \max_{\theta\in\NN} \tfrac{(\beta^{\theta}_f)^2}{8} \| \EE_\theta(P^I) \|\delta^2 + \tfrac{\beta^{\theta}_f}{2} \| \Gamma_{\theta} \| \| \EE_\theta(P^I)  \|\delta < 1,
\end{align*}
and $\alpha$ so that
\begin{align*}
 \alpha \geq \max_{\theta\in\NN} \tfrac{\frac{(\beta^{\theta}_f)^2}{8} \| \EE_\theta(P^\beta) \| \delta^2  + \frac{\beta^{\theta}_f}{2} \| \Gamma_{\theta} \| \| \EE_\theta(P^\beta)  \|\delta + \frac{\beta^{\theta}_f}{2} \|\EE_\theta(p) \|}
 {1-\frac{(\beta^{\theta}_f)^2}{8} \| \EE_i(P^I) \|\delta^2 + \frac{\beta^{\theta}_f}{2} \| \Gamma_{\theta} \| \| \EE_\theta(P^I)  \|\delta}.
\end{align*}
Now for $x\in\mathcal{B}_{\delta}$ and $\alpha$ as above 
we have
$\Delta\mathcal{L}V_f(x,\theta) \leq \tfrac{\alpha}{2}\|x\|^2$, 
and since $\mathcal{L}V_f(x,\theta) = -\ell_q(x,\theta) 
+ \Delta \mathcal{L}V_f(x,\theta)$ we have
\begin{align*}
    \mathcal{L} V_f(x,\theta) \le -\ell_q(x,\theta) 
  + \tfrac{\alpha}{2}\|x\|^2,
\end{align*}
and employing Lemma~\ref{lemma:l_q} we obtain
\begin{align}
 \mathcal{L} V_f(x,\theta)\le-\bar{\ell}(x,\theta).
\end{align}

If Assumption~\ref{assum:lower-quad-bound} holds all 
assumptions of Theorem~\ref{thm:mss_of_terminal_set_MPC} are 
fulfilled and the controlled system is locally mean square 
stable.}~\hfill{$\Box$}
\end{pf} % ------ Proof `thm:lin-stability` ------

\subsection{Computation of $X^f$}
We shall demonstrate one possible way of finding $X^f$ such 
that the requirements of Theorem~\ref{thm:lin-stability} 
are satisfied. Take $X^f=\{X^f_i\}_{i\in\NN}$ to be ellipsoidal 
of the form $X^f_i=\{x: x^\top P_i x \leq 1\}$.
By Assumption~\ref{assum:beta_smooth_system}, there are 
constants $\gamma_i> 0$, $i\in\NN$, so that 
\begin{align}
 x_{k+1} &= 
   A_{\theta_k} x_k 
 + B_{\theta_k} \kappa_f(x_k, \theta_k) 
 + d_{k, \theta_k},
\end{align}
with $\| d_{k,i} \|^2 \le \gamma_i  x_k^\top P^f_i x_k$ where 
$d_{k,i}=e(x_k,i)$ is the linearization error.
For $X^f$ to be UPI for the $\kappa_f$-controlled system 
it must satisfy
\begin{subequations}
\begin{align}
&\max_{j \in \mathcal{C}(i)} \{x_{k+1}^\top P_j x_{k+1} \}\le x_k^\top P_i x_k, \,\forall i \in \NN\notag\\
\Leftrightarrow\,&
   \begin{bmatrix}x_k\\d_{k,i}\end{bmatrix}^\top
   \begin{bmatrix}
    P_i - \Gamma_i^\top P_j \Gamma_i & -\Gamma_i^\top P_j\\
    -P_j\Gamma_i & -P_j
    \end{bmatrix}
    \begin{bmatrix}x_k\\d_{k,i}\end{bmatrix}
    \geq 0,\label{eq:desired_decrease}
\end{align}
for all $ j \in \mathcal{C}(i)$ and $i\in\NN$
whenever $d_{k,i}^\top d_{k,i} \leq \gamma_i x_k^\top P^f_i x_k$, or, for $i\in\NN$
\begin{align}\label{eq:lmi_premise}
 \begin{bmatrix}x_k\\d_{k,i}\end{bmatrix}^\top
 \begin{bmatrix}\gamma_i P_i^f\\& -I\end{bmatrix}
 \begin{bmatrix}x_k\\d_{k,i}\end{bmatrix}\ge 0.
\end{align}
\end{subequations}
Using the S\,-lemma,~\eqref{eq:lmi_premise} implies~\eqref{eq:desired_decrease} 
so long as 
\begin{align}\label{eq:x01}
\begin{bmatrix}
P_i - \Gamma_i^\top P^f_j \Gamma_i & -\Gamma_i^\top P^f_j\\
-P_j\Gamma_i & -P^f_j
\end{bmatrix}
- \tau \begin{bmatrix}\gamma_i P^f_i\\& -I\end{bmatrix} \succcurlyeq 0
\end{align}
for some $\tau\geq 0$ and for all $i\in\NN$ and $j\in\mathcal{C}(i)$.
By rearranging the terms in the two matrices, equation~\eqref{eq:x01} can be equivalently written as
\begin{align}\label{eq:x02}
\begin{bmatrix}\tau \gamma_i P^f_i + \Gamma_i ^\top  P^f_j \Gamma_i & \Gamma_i^\top P^f_j\\ * & P^f_j\end{bmatrix}
\preccurlyeq
 \begin{bmatrix}P^f_i \\  & \tau I \end{bmatrix}.
\end{align}
The left hand side of~\eqref{eq:x02} is equal to
\begin{align*}
 \begin{bmatrix}P^f_i\\ P^f_j\Gamma_i & P^f_j\end{bmatrix}^\top 
 \begin{bmatrix}\tfrac{1}{\gamma_i \tau}P^f_i\\ & P^f_j\end{bmatrix}^{-1}
 \begin{bmatrix}P^f_i\\ P^f_j\Gamma_i & P^f_j\end{bmatrix}
\end{align*}
Using the Schur complement we get 
\begin{align}\label{eq:x03}
 \begin{bmatrix}
  P^f_i & 0 & P^f_i & \Gamma_i^\top P^f_j\\
  0   & \tau I & 0 & P^f_j\\
 * & * & \tfrac{1}{\gamma_i \tau}P^f_i&0\\
 * & * & 0 & P^f_j 
 \end{bmatrix} \succcurlyeq 0.
\end{align}
Introducing the variables $P^f_i = Z_i^{-1}$ and 
$K_i = Y_iZ_i^{-1}$,~\eqref{eq:x02} 
is equivalent to the matrix inequality
%
%
% Final LMI --- #1
%
\begin{align}\label{eq:x04}
 \begin{bmatrix}
  Z_i & 0 & \tau Z_i & Z_i A_i^\top + Y_i^\top B_i^\top\\
  0   & \tau I & 0 & I\\
 * & * & \tau \gamma_i^{-1}Z_i & 0\\
 * & * & 0 & Z_j 
 \end{bmatrix}  \succcurlyeq 0.
\end{align}
As required by Theorem~\ref{thm:lin-stability}, $X_i^f$ must 
be in $\mathcal{B}_\delta$.
This is equivalently written as
%
%
% Final LMI --- #2
%
\begin{align}\label{eq:x05}
 \begin{bmatrix}
  \delta I & P_i\\
  P_i & I
 \end{bmatrix} \succcurlyeq 0.
\end{align}
We then choose $P^f_i$ so as to satisfy~\eqref{eq:x04} and~\eqref{eq:x05}
for all $i\in\NN$ and $j\in\mathcal{C}(j)$.
Note that~\eqref{eq:x04} is a bilinear matrix inequality (BMI) with 
unknowns $Z_i$, $Y_i$ and $\tau$, but the bilinearity is only
because of the term $\tau Z_i$.
Although BMIs are more difficult to solve compared to LMIs, 
in this case since $\tau$ is a scalar,~\eqref{eq:x04} can be 
solved with a simple line search method with respect to $\tau$.

\section{Conclusions}
This paper offers a theoretical framework for the control 
of Markovian switching systems using EMPC. 
We first studied a formulation with mode-dependent 
optimal steady states and terminal equality constraints for which 
we provided an upper bound  on the expected asymptotic average cost 
(Theorem~\ref{thm:asymptotic-perf}). 
We then studied an EMPC formulation with mode-dependent terminal 
region constraints and we provided design guidelines based on 
the system linearization assuming that the system dynamics 
and the stage cost function are $\beta$-smooth which are rather weak 
assumptions
(Theorem~\ref{thm:lin-stability}).

% REFERENCES
\bibliography{ifacconf}

\end{document}